\documentclass[12pt]{amsart}

\usepackage[dvips]{graphics}

\usepackage{epsfig}

\usepackage{amsfonts}
\usepackage{amsmath, amsthm, amssymb, ulem, amscd} 

\def\crn#1#2{{\vcenter{\vbox{
        \hbox{\kern#2pt \vrule width.#2pt height#1pt
           }
          \hrule height.#2pt}}}}
\def\intprod{\mathchoice\crn54\crn54\crn{3.75}3\crn{2.5}2}
\def\into{\mathbin{\intprod}}

\newcommand{\pa}{\partial}
\newcommand{\Up}{\Upsilon}

\newcommand{\mut}{{\tilde{\mu}}}
\newcommand{\cG}{{\mathcal G}}
\newcommand{\cA}{{\mathcal A}}
\newcommand{\cC}{{\mathcal C}}

\newcommand{\cM}{{\mathcal M}}
\newcommand{\cGt}{{\tilde {\mathcal G}}}

\newcommand{\cO}{{\mathcal O}}
\newcommand{\cOt}{\widetilde{\mathcal O}}

\newcommand{\gt}{{\tilde g}}

\newcommand{\Rt}{{\tilde R}}

\newcommand{\R}{{\mathbb R}}

\newcommand{\C}{{\mathbb C}}

\newcommand{\nt}{\tilde \nabla}

\newcommand{\Ric}{\operatorname{Ric}}
\newcommand{\contr}{\operatorname{contr}}
\newcommand{\tf}{\operatorname{tf}}
\renewcommand{\Re}{\mathop{\rm Re}\nolimits}

\renewcommand{\tilde}{\widetilde} 
\newcommand{\wh}{\widehat} 

\theoremstyle{plain}
\newtheorem{theorem}{Theorem}[section]

\newtheorem{proposition}[theorem]{Proposition}

\theoremstyle{definition}

\newtheorem{definition}[theorem]{Definition}

\theoremstyle{remark}

\numberwithin{equation}{section}

\title{Inhomogeneous Ambient Metrics} 

\author{C. Robin Graham}
\address{Department of Mathematics, University of Washington,
Box 354350\\
Seattle, WA 98195-4350}
\email{robin@math.washington.edu}

\author{Kengo Hirachi}
\address{Graduate School of Mathematical Sciences,
University of Tokyo\\
3-8-1 Komaba, Meguro\\
Tokyo 153-8914, Japan}
\email{hirachi@ms.u-tokyo.ac.jp}

\begin{document}

\maketitle

\thispagestyle{empty}

\renewcommand{\thefootnote}{}
\footnotetext{This research was partially supported by NSF grant \# DMS 
  0505701 and Grants-in-Aid for Scientific Research, JSPS.}   
\renewcommand{\thefootnote}{1}

\section{Introduction}\label{intro}

The ambient metric, introduced in \cite{FG1}, has proven to be an important
object in  
conformal geometry.  To a manifold $M$ of dimension $n$ with
a conformal 
class of metrics $[g]$ of signature $(p,q)$ it associates a formal
expansion of a metric 
$\gt$ of signature $(p+1,q+1)$ on a space $\cGt$ of dimension $n+2$.  
This generalizes the realization of the conformal sphere $S^n$ as the space
of null lines for a quadratic form of signature $(n+1,1)$, with associated  
Minkowski metric $\gt$ on $\R^{n+2}$.  The ambient space $\cGt$ carries a
family of 
dilations with respect to which $\gt$ is homogeneous of 
degree 2.  The other conditions determining $\gt$ are that it
be Ricci-flat and satisfy an initial condition specified by the   
conformal class $[g]$.  

The ambient metric behaves differently in even and odd dimensions,
reflecting an underlying distinction in the structure of the space of jets
of metrics modulo 
conformal rescaling.  When $n$ is odd, the equations
determining $\gt$ have a smooth infinite-order formal power 
series solution, uniquely determined modulo diffeomorphism.  But when $n$
is even and $\geq 4$, there is a obstruction at order $n/2$ to the
existence of a smooth formal solution, which is realized as a conformally  
invariant natural tensor called the ambient obstruction
tensor.  It is possible to continue the expansion of $\gt$ to
higher orders by 
including log terms (\cite{K}), but this destroys the differentiability of
the solution so is problematic for applications requiring 
higher order differentiation.

In this article, we describe a modification to the form of the ambient
metric in even dimensions which enables us to obtain invariantly defined,
smooth infinite order ``ambient metrics''.  We introduce what we call
inhomogeneous ambient metrics, which are formally Ricci-flat and have 
an asymptotic 
expansion involving the logarithm of a defining function for the
initial surface which is homogeneous of degree 2.  Such metrics are
themselves no longer homogeneous, and of course are not
smooth.  However, we are able to define the smooth part of such a
metric in an invariant way, and the smooth part is homogeneous 
and of course smooth and can 
be used in applications just as the ambient metric itself is used in odd
dimensions.  A significant difference, though, is that an inhomogeneous
ambient metric is no longer uniquely
determined up to diffeomorphism by the conformal class $[g]$ on $M$: there 
is a family of inhomogeneous ambient metrics, and therefore of their smooth
parts, 
parametrized by the choice of an arbitrary trace-free symmetric 2-tensor
field on $M$ which we call the ambiguity tensor.  

There are two main motivations for the introduction of terms involving the
logarithm of a defining function homogeneous of degree 2 in the expansion
of an ambient metric.  The simplest is the observation that in flat 
space, the null cone of the Minkowski metric has an invariant 
defining function homogeneous of degree 2, namely the quadratic form $Q$ 
itself, but no invariant defining function homogeneous of degree 0.  
In this flat space setting, the theory of invariants of conformal densities 
initiated in \cite{EG} can be extended to all orders by the introduction of 
log $|Q|$ terms in the expansion of the ambient harmonic extension of a
density.  We will report on this work on invariants of densities elsewhere.      

Another motivation is the construction of the  
potential of the inhomogeneous CR ambient metric in \cite{H}, which also
involves the log 
of a defining function homogeneous of degree 2 and an invariantly defined
smooth part.  The existence of this construction in the CR case was
compelling evidence  
that there should be a conformal analogue.  Nonetheless, despite having the
CR construction as a    
guide, we did not find the conformal construction to be straightforward.
The relation between the constructions is still not completely
clear; further discussion of this issue is contained in \S\ref{inhomo}.   

We also indicate how inhomogeneous ambient  
metrics can be used to complete the description of scalar invariants of   
conformal structures in even dimensions.  One can form scalar invariants as 
Weyl invariants, defined to be linear combinations of complete 
contractions of covariant derivatives of the curvature tensor of the
smooth part of an inhomogeneous ambient metric, and also its volume
form and a modified volume form in the case of odd invariants.  Such
invariants generally depend  
on the choice of ambiguity.  Nonetheless, in dimensions $n \equiv 2 \mod
4$, all scalar conformal invariants arise as Weyl invariants which are
independent of the ambiguity.  If $n \equiv 0 \mod 4$, 
there are exceptional odd invariants which are not of this form.  In this
case, following \cite{BG}, we provide a construction of a finite set of 
basic exceptional invariants, which, together with ambiguity-independent 
Weyl invariants, span all scalar conformal invariants.  
The study of scalar invariants was initiated in 
\cite{F2} in the CR case.  Fundamental algebraic results were derived in 
in \cite{BEGr}, resulting in the full 
description of scalar invariants in odd dimensional conformal geometry and
below the 
order of the obstruction in CR and even dimensional conformal geometry.  
The  
completion of the description in the CR case was given in \cite{H}. 
A different approach to the study of conformal invariants is developed in
\cite{Go} using tractor calculus.  
In \cite{A1}, \cite{A2}, Alexakis has recently derived a theory handling
invariants of a density coupled with a conformal structure, below the order  
of the obstruction of the ambient metric in even dimensions, and below the 
order of the obstruction of harmonic extension of the density.  

In \S\ref{homo}, we review the construction of smooth homogeneous ambient
metrics.  Inhomogeneous ambient  
metrics are introduced in \S\ref{inhomo} and the main theorem asserting the 
existence and uniqueness of inhomogeneous ambient metrics with prescribed 
ambiguity tensor is formulated.  We describe in some detail the
transformation law for the ambiguity tensor under conformal change.  We  
state a necessary and sufficient condition for the asymptotic expansion of 
an inhomogeneous ambient metric to have a simpler form, and show by 
explicit identification of the obstruction tensor 
that Fefferman metrics of nondegenerate integrable CR manifolds 
always satisfy this condition.  We also briefly discuss the Poincar\'e
metrics associated to inhomogeneous ambient metrics.  In
\S\ref{invariants}, we describe the results concerning scalar invariants
and give some examples of invariants.

\section{Smooth homogeneous ambient metrics}\label{homo} 

In this section we review the usual ambient metric construction in odd
dimensions and up to the obstruction in even dimensions.  Let $[g]$ be a 
conformal class of metrics of signature $(p,q)$ on a (paracompact) manifold
$M$ of 
dimension $n\geq 3$.  Here $g$ is a smooth metric of signature $(p,q)$ 
on $M$ and $[g]$ consists of all metrics of the form $\Omega^2 g$ with  
$0<\Omega \in C^{\infty}(M)$.  The metric bundle 
$\cG\subset \bigodot^2T^*M$ of $[g]$ is the set of all pairs  
$(x,\underline{g})$, where $x\in M$ and 
$\underline{g}\in \bigodot^2T_x^*M$ is of the form $\underline{g}=t^2g(x)$
for some $t>0$.  The fiber variable $t$ on $\cG$ so defined is  
associated to the choice of 
metric $g$ and provides an identification $\cG\cong \R_+\times M$.   
Sections of $\cG$ are precisely metrics in the conformal class.   
There is a tautological symmetric $2$-tensor $g_0$ on $\cG$ defined by 
$g_0(X,Y)=\underline{g}(\pi_*X,\pi_*Y)$, where $\pi:\cG\rightarrow M$ is the
projection 
and $X$, $Y$ are tangent vectors to $\cG$ at $(x,\underline{g})\in \cG$.
The family of dilations $\delta_s:\cG\rightarrow \cG$ defined by  
$\delta_s(x,\underline{g})= (x,s^2\underline{g})$ defines an $\R_+$ action
on $\cG$, and one has  
$\delta^*g_0 = s^2g_0$.  
We denote by $T= \frac{d}{ds} \delta_s |_{s=1}$ the vector field on $\cG$ which is
the infinitesimal generator of the dilations $\delta_s$.

The ambient space is $\cGt= \cG\times \R$;  
the coordinate in the $\R$ factor is typically written $\rho$.
The dilations $\delta_s$ extend to $\cGt$ acting on the $\cG$ factor
and we denote also by $T$ the infinitesimal generator of the $\delta_s$ on 
$\cGt$. 
We embed $\cG$ into $\cGt$ by $\iota: z\rightarrow (z,0)$ for $z\in
\cG$, and we identify $\cG$ with its image under $\iota$.  As described
above, a choice of representative metric $g$ induces 
an identification $\cG \cong \R_+\times M$, so also an identification
$\cGt\cong \R_+\times M\times \R$.  We use capital Latin indices to 
label objects on $\cGt$, '$0$' indices to label the $\R_+$ 
factor, lower case Latin indices to label the $M$ factor, and '$\infty$' 
indices for the $\R$ factor.   

In the even-dimensional case, different components of the ambient metric
are determined to different orders.  
If $S_{IJ}$ is a smooth symmetric 2-tensor field on an   
open neighborhood of $\cG$ in $\cGt$, then 
for $m \geq 0$ we write $S_{IJ} = O^+_{IJ}(\rho^m)$ if:   
\begin{enumerate}
\item[(i)]
$S_{IJ} = O(\rho^m)$; and
\item[(ii)]
For each point $z\in \cG$, the symmetric 2-tensor $(\iota^*(\rho^{-m}S))(z)$  
is of the form 
$\pi^*s$ for some symmetric 2-tensor $s$ at $x=\pi(z)\in M$ satisfying 
$\operatorname{tr}_{g_x}s = 0$.  The symmetric 2-tensor $s$ is allowed to
depend on $z$, not just on $x$.  
\end{enumerate}
In terms of components relative to a choice of representative metric $g$,
$S_{IJ}= O^+_{IJ}(\rho^m)$ if and only if  
all components satisfy $S_{IJ}= O(\rho^m)$ and if in addition one has that 
$S_{00}$, $S_{0i}$ and $g^{ij}S_{ij}$ are $O(\rho^{m+1})$.  The condition
$S_{IJ}=O^+_{IJ}(\rho^m)$ is easily seen to be preserved by diffeomorphisms  
on a neighborhood of $\cG$ 
which restrict to the identity on $\cG$.

We consider metrics $\gt$ of signature $(p+1,q+1)$ defined on a
neighborhood of $\cG$ in $\cGt$.  We will assume that the neighborhood is
homogeneous, i.e., it is invariant under the dilations $\delta_s$ for  
$s>0$.  A smooth diffeomorphism defined on such a neighborhood will be said
to be homogeneous if it commutes with the $\delta_s$.  The metric $\gt$
(or more generally a symmetric 2-tensor field) will be said to be
homogeneous (of degree 2) if it satisfies $\delta_s^*\gt=s^2\gt$.   

The main existence and uniqueness result for smooth ambient metrics is the
following.
\begin{theorem}\label{beforeobs}
Let $[g]$ be a conformal class on a manifold $M$.  If $n$ is odd, then
there is a smooth metric $\gt$ on a homogeneous
neighborhood of $\cG$ in $\cGt$, uniquely determined  
up to a homogeneous diffeomorphism of a neighborhood of $\cG$
in $\cGt$ which restricts to the
identity on $\cG$, and up to a homogeneous term vanishing to infinite order  
along $\cG$, by the requirements:
\begin{enumerate}
\item
$\delta_s^*\gt = s^2\gt$
\item
$\iota^*\gt = g_0$
\item
$\Ric(\gt) = 0$ to infinite order along $\cG$.
\end{enumerate}
If $n\geq 4$ is even, the same statement holds except that (3) is replaced
by 
$\Ric(\gt) = O^+_{IJ}(\rho^{n/2-1})$, and $\gt$ is uniquely determined 
up to a homogeneous diffeomorphism and up to a smooth homogeneous term
which is $O^+_{IJ}(\rho^{n/2})$.      
\end{theorem}

\noindent
A metric $\gt$ satisfying the conditions of Theorem~\ref{beforeobs} is
called a smooth ambient metric for $(M,[g])$.

The diffeomorphism indeterminacy in $\gt$ can be fixed by the choice of a 
metric $g$ in the conformal class.  As described above, the choice of such
a metric $g$ determines an identification $\cGt\cong \R_+\times M \times
\R$.   
\begin{definition}\label{normalformdef}
A smooth metric $\gt$ on $\cGt$ satisfying (1) and (2) 
above is said to be in normal form relative to $g$ if in the identification
$\cGt = \R_+\times M\times \R$ induced by $g$,
\begin{enumerate}
\item[(a)]
$\gt = 2t\,dt\,d\rho + g_0$ at $\rho = 0$, and 
\item[(b)]
The lines 
$\rho\rightarrow (t,x,\rho)$ are geodesics for $\gt$ for each choice of 
$(t,x)\in \R_+\times M$.  
\end{enumerate}
\end{definition}

As in the construction of Gaussian coordinates relative to a  
hypersurface, it follows by straightening out geodesics that if $\gt$ is
any smooth metric satisfying (1) and (2) and $g$ is a  
representative metric in the conformal class, then there is a unique
homogeneous diffeomorphism $\Phi$ defined in a homogeneous neighborhood of
$\cG$ and which restricts to the identity on $\cG$, such that $\Phi^*\gt$   
is in normal form relative to $g$.  It follows that in
Theorem~\ref{beforeobs}, the metric $\gt$ can be chosen to be in 
normal form relative to $g$, and it is then uniquely determined up to 
$O(\rho^{\infty})$ for $n$ odd, and up to $O^+_{IJ}(\rho^{n/2})$ for $n$
even.  

Theorem~\ref{beforeobs} is proved by a formal power series analysis of the 
equation $\Ric(\gt)=0$ for metrics $\gt$ in normal form relative to some
metric $g$ in the conformal class.  In carrying out this analysis, one
finds that the solution has an additional property.  

\begin{definition}
A metric $\gt$ on a homogeneous neighborhood $\mathcal{U}$ of $\cG$ in
$\cGt$ is said to be {\it straight} if for each $p\in \mathcal{U}$, the
dilation orbit $s\rightarrow \delta_sp$ is a geodesic for $\gt$.  
\end{definition}

\begin{proposition}\label{straight}
In Theorem~\ref{beforeobs}, the metric $\gt$ can be chosen to be straight.
\end{proposition}

\noindent
Since the condition that $\gt$ be straight is invariant under smooth
homogeneous diffeomorphisms, it follows that any smooth ambient metric
$\gt$ agrees with a straight metric to infinite order 
when $n$ is odd and mod $O^+_{IJ}(\rho^{n/2})$ when $n$ is even.  That is,
any smooth ambient metric is straight to the order that it is determined.

When $n$ is even, in general there is an obstruction to the existence of 
a smooth ambient metric solving $\Ric(\gt)=O(\rho^{n/2})$.
Let $\gt$ be a smooth ambient metric.
We denote by $r_\#$ the function $r_\#=\|T\|^2 =\gt(T,T)$.  Since $g_0(T,T)=0$,
it follows that $r_\#=0$ on $\cG$, and it turns out (as a consequence of the
Einstein condition or the straightness condition) that  
$dr_\#\neq 0$ on $\cG$.  Thus $r_\#$ is a defining function for $\cG\subset \cGt$
invariantly associated to $\gt$ which is homogeneous of degree 2:
$\delta_s^*r_\#=s^2r_\#$.  
Now since $\Ric(\gt)=O^+_{IJ}(\rho^{n/2-1})$, the quantity
$(r_\#^{1-n/2}\Ric{\gt})|_{T\cG}$ defines a tensor field on $\cG$, 
homogeneous of degree $2-n$, which annihilates $T$.  It therefore defines a 
symmetric 2-tensor-density on $M$ of weight $2-n$, which is trace-free.
If $g$ is a metric in the conformal class, evaluating this tensor-density
at the image of $g$ viewed as a section of $\cG$ defines a 2-tensor on $M$
which we denote by 
$(r_\#^{1-n/2} \Ric{\gt})|_g$.  The obstruction tensor of the metric $g$ is
defined to be  
\begin{equation}\label{obsdef}
\cO = c_n (r_\#^{1-n/2} \Ric{\gt})|_g, 
\qquad c_n = (-1)^{n/2-1} \frac{2^{n-2} (n/2-1)!^2}{n-2}.
\end{equation}
We then have:

\begin{proposition}\label{obstruction}
Let $n\geq 4$ be even.  The obstruction tensor $\cO_{ij}$ of $g$ is
independent of 
the choice of smooth ambient metric $\gt$ and has the following properties:    
\begin{enumerate}
\item
$\cO$ is a natural tensor invariant of the metric $g$;
i.e. in local coordinates the components of $\cO$ are given by
universal polynomials in the components of $g$, $g^{-1}$ and the curvature 
tensor of $g$ and its covariant derivatives.
The expression for $\cO_{ij}$  takes the form 
\begin{equation}\label{Oform}
\begin{split}
\cO_{ij} &= \Delta^{n/2-2}\left (P_{ij},_k{}^k-P_k{}^k,_{ij}\right )
+ lots\\
&=(3-n)^{-1}\Delta^{n/2-2}W_{kijl},{}^{kl} + lots,
\end{split}
\end{equation}
where  
$\Delta = \nabla^i\nabla_i$, 
$$
P_{ij}=\frac{1}{n-2}\left(R_{ij}-\frac{R}{2(n-1)}g_{ij}\right),
$$
$W_{ijkl}$ is the Weyl tensor,  
and lots denotes quadratic and higher terms in curvature involving fewer 
derivatives.  
\item
One has
$$
\cO_i{}^i=0 \qquad\qquad\qquad \cO_{ij},{}^j = 0.
$$
\item
$\cO_{ij}$ is conformally invariant of weight $2-n$; i.e. if $0<\Omega \in 
  C^{\infty}(M)$ and $\wh{g}_{ij} =
\Omega^2 g_{ij}$, then $\wh{\cO}_{ij} = \Omega^{2-n}\cO_{ij}$.   
\item
If $g_{ij}$ is conformal to an Einstein metric, then $\cO_{ij}=0$.  
\item 
If $n=4$, then $\cO_{ij} = B_{ij}$ is the classical Bach tensor.
\end{enumerate}
\end{proposition}

\noindent
Here the Bach tensor is defined in any dimension $n\geq 3$ by
$$
B_{ij} = C_{ijk},{}^k - P^{kl}W_{kijl},
$$
where $C_{ijk}$ is the Cotton tensor
$$
C_{ijk} = P_{ij},_k-P_{ik},_j.
$$

Clearly, if the obstruction tensor is nonzero, then it is impossible to
find a smooth ambient metric $\gt$ solving $\Ric(\gt) = O(\rho^{n/2})$.

\section{Inhomogeneous ambient metrics}\label{inhomo}
Restrict now to the case $n$ even.  In order to find ambient metrics beyond
order $n/2$, we broaden the class of allowable metrics $\gt$.  The $\gt$
that we consider are neither homogeneous nor smooth.  However, $\gt$ will
always be required to satisfy the initial condition
\begin{equation}\label{initial}
\iota^*\gt = g_0.
\end{equation}

Let $r$ be a smooth defining function for 
$\cG\subset \cGt$ homogeneous of degree 0 (one could for example choose
$r=\rho$) and 
let $r_\#$ be a smooth defining function homogeneous of degree 2 (for
example, $r_\# = 2\rho t^2$).  Denote by $\cM$ the space of formal asymptotic
expansions along $\cG$ of metrics on $\cGt$ of signature $(p+1,q+1)$ of the
form    
\begin{equation}\label{metricexpansion}
\gt \sim \gt^{(0)}+ \sum_{N\geq 1}\gt^{(N)}r (r^{n/2-1}\log |r_\#|)^N    
\end{equation}
where each $\gt^{(N)}$, $N\geq 0$, is a smooth symmetric 2-tensor field on 
$\cGt$ 
satisfying $\delta_s^* \gt^{(N)}=s^2\gt^{(N)}$, such that the initial 
condition \eqref{initial} holds.  Observe that $\gt \in\cM$ is the
expansion of a smooth metric on $\cGt$ if
and only if $\gt^{(N)}=0$  to infinite order for $N\geq 1$.
We will say in this case that $\gt$ is smooth.  Note that if $\gt$ is
smooth, then $\gt$ agrees to infinite order with a homogeneous metric; we
will also say that $\gt$ is homogeneous. 
Denote by $\cA$ the space of formal asymptotic expansions of  
scalar functions $f$ on $\cGt$ of the form  
$$
f\sim f^{(0)}+ \sum_{N\geq 1} f^{(N)}r (r^{n/2-1}\log |r_\#|)^N, 
$$
where each $f^{(N)}$ is a smooth function on $\cGt$ homogeneous of
degree $0$.  Observe that these spaces of
asymptotic expansions are independent of the choice of defining
functions $r$, $r_\#$.  
If $\Phi$ is a formal smooth, homogeneous diffeomorphism of $\cGt$
satisfying $\Phi|_\cG = Id$, then pullback by $\Phi$ preserves   
$\cM$ and $\cA$.  
As already initiated above, in the following we will refer to  
such asymptotic expansions for metrics, functions, and formal smooth
diffeomorphisms as if 
they were actually defined in a neighborhood of $\cG$,  
and when we say that such an object satisfies a certain condition,
we will mean that the condition holds formally to infinite order along
$\cG$.  For example, we will say that $\gt$ is straight if  
the ordinary differential equations which state that the orbits 
$s\rightarrow \delta_sp$ are geodesics hold formally to infinite order 
along $\cG$.  

For inhomogeneous ambient metrics, we impose the straightness condition at
the outset.   

\begin{definition}\label{ambmetric}
An inhomogeneous ambient metric for $(M,[g])$ is a straight metric $\gt\in
\cM$ satisfying $\Ric(\gt)=0$.     
\end{definition}

The straightness 
condition is crucial in the inhomogeneous case because of the following
proposition. 

\begin{proposition}\label{Tsmooth}
Let $\gt\in \cM$ be straight.  Then $T\into \gt$ is 
smooth and $\gt(T,T)$ is a smooth defining function for $\cG$ homogeneous
of degree 2.
\end{proposition}

Observe that for $\gt\in \cM$, $T\into \gt$ and $\gt(T,T)$ will in   
general have asymptotic expansions involving $\log |r_\#|$, so will be
neither smooth nor homogeneous.  Proposition~\ref{Tsmooth} asserts that the 
requirement that $\gt$ be straight has as a 
consequence that no log terms occur in the expansions of these quantities.  
The proof of Proposition~\ref{Tsmooth} is a straightforward analysis of the
geodesic equations for the dilation orbits.

If $\gt\in \cM$ is straight, then $\gt(T,T)$ is a 
canonically determined smooth defining function for $\cG$ 
homogeneous of degree 2.  We may therefore take $r_\#=\gt(T,T)$ in
\eqref{metricexpansion}.  Having fixed $r_\#$, the smooth symmetric
2-tensor fields $\gt^{(0)}$ and $r^{(n/2-1)N+1}\gt^{(N)}$ for $N\geq 1$ are 
then uniquely determined by $\gt$    
independently of any choices.  In particular, $\gt^{(0)}$ is an invariantly  
determined smooth part of $\gt$, which is itself a smooth metric in $\cM$.  
If $\Phi$ is a smooth homogeneous diffeomorphism satisfying $\Phi |_\cG =  
Id$ and $\gt\in\cM$ is straight, then $(\Phi^*\gt)^{(0)}=
\Phi^*(\gt^{(0)})$.  

Next we formulate the notion of normal form for straight metrics in $\cM$.   
\begin{definition}
Let $g$ be a metric in the conformal class and let $\gt\in\cM$ be
straight.  Then $\gt$ is said to be in normal form relative to $g$ if its
smooth part is in normal form relative to $g$.   
\end{definition}
\noindent
It is clear from the existence and uniqueness of a diffeomorphism putting a
smooth metric into normal form that if $\gt\in\cM$ is straight, then there
is a unique  
smooth homogeneous diffeomorphism $\Phi$ such that $\Phi |_\cG = Id$ and
such that $\Phi^*\gt$ is in normal form relative to $g$.    

The following proposition makes explicit the normal form condition.  Its
proof is an analysis of the geodesic equations for the straightness and 
normal form conditions using the initial condition (a) of 
Definition~\ref{normalformdef}.  
\begin{proposition}\label{normalformprop}
Let $g$ be a representative for the conformal structure.  A straight metric 
$\gt \in \cM$ is in normal form relative to $g$ if and only if in the  
identification $\cGt = \R_+\times M\times \R$ induced by $g$, $\gt$ takes
the form: 
\begin{equation}\label{formgt}
\gt_{IJ} =  
\begin{pmatrix}
2\rho&0&t\\
0&t^2g_{ij}&t^2g_{i\infty}\\
t&t^2g_{j\infty}&t^2g_{\infty\infty}
\end{pmatrix}
\qquad
\text{ with $g_{ij}$, $g_{i\infty}$, $g_{\infty\infty}\in \cA$}, 
\end{equation}
where $g_{ij}|_{\rho = 0}$ is the given metric $g$ on $M$, 
and where the expansions for 
the components $g_{j\infty}$ and $g_{\infty\infty}$ have zero smooth part 
when expanded using $r_\# = 2\rho t^2$.  
That is, for $I= i$, $\infty$ we have  
\begin{equation}\label{normalvanishing}
g_{I\infty}\sim \sum_{N\geq 1}g_{I\infty}^{(N)}
\rho (\rho^{n/2-1}\log |2\rho t^2|)^N 
\end{equation}
where the $g_{I\infty}^{(N)}$ are smooth and homogeneous of degree 0.    
\end{proposition}

It is a consequence of Proposition~\ref{normalformprop} that 
if $\gt\in\cM$ is straight, then its smooth part is also straight.

The next theorem is our main theorem concerning the existence and
uniqueness of inhomogeneous ambient metrics.  As described in \S\ref{homo}, 
in odd dimensions, given a representative metric $g$ in the conformal 
class, there is to infinite order a unique smooth ambient metric in normal
form relative to $g$.  For $n$ even, inhomogeneous ambient metrics in
normal form are no longer uniquely determined by a representative for the
conformal class.  One has exactly the freedom of a smooth trace-free
symmetric 2-tensor field on $M$, which we call the ambiguity tensor.  

\begin{theorem}\label{main}
Let $(M^n,[g])$ be a manifold with a conformal structure, $n\geq 4$ even.    
Up to pull-back by a smooth homogeneous diffeomorphism which
restricts to the identity on $\cG$, the 
inhomogeneous ambient metrics for $(M,[g])$ are parametrized by the choice 
of an arbitrary trace-free symmetric 2-tensor field $A_{ij}$ on $M$.   
Precisely, for each representative metric $g$ and choice of ambiguity 
tensor $A_{ij}\in \Gamma(\bigodot^2_0T^*M)$, there is a unique  
inhomogeneous ambient metric $\gt$ in normal form relative to $g$ such that  
\begin{equation}\label{ambig}
\tf{\left( (\pa_\rho)^{n/2}g_{ij}^{(0)}|_{\rho =0}\right)} 
=  A_{ij}.    
\end{equation}
Here we have written $\gt$ in
the form \eqref{formgt}, $g^{(0)}_{ij}$ denotes the
smooth part of $g_{ij}$, and $\tf$ the trace-free part.

There is a natural pseudo-Riemannian invariant  1-form $D_i$ so that the  
metric $\gt$ determined by initial metric $g$ and ambiguity $A_{ij}$ is 
smooth if and only if $\cO_{ij}=0$ and $A_{ij},^j=D_i$, where   
$\cO_{ij}$ denotes the obstruction tensor of $g$.    
\end{theorem}

The $1$-form $D_i$ is defined as follows.  A straight smooth ambient 
metric in normal form takes the form 
\begin{equation}\label{smoothform}
\gt_{IJ}=
\begin{pmatrix}
2\rho&0&t\\
0&t^2g_{ij}&0\\
t&0&0
\end{pmatrix}
\end{equation}
where $g_{ij}=g_{ij}(x,\rho)$ is a smooth 1-parameter family of metrics on
$M$.  The derivatives $\pa_\rho^mg_{ij}$ for $m\leq n/2-1$ and
$g^{ij}\pa_\rho^{n/2}g_{ij}$ are determined at $\rho =0$ and are natural
invariants of the initial metric $g_{ij}(x,0)$.  Fix the indeterminacy
of $g_{ij}(x,\rho)$ to higher orders by fixing $g_{ij}(x,\rho)$ to be the  
Taylor polynomial of degree $n/2$ whose $(n/2)^{nd}$ derivative is 
pure trace:
\begin{equation}\label{expand}
g_{ij}(x,\rho)=\sum_{m=0}^{n/2-1}\frac{1}{m!}\pa_\rho^mg_{ij}\rho^m   
+\frac{1}{n(n/2)!}(g^{kl}\pa_\rho^{n/2}g_{kl})g_{ij}\rho^{n/2},
\end{equation}
where on the right hand side, $g_{ij}$ and its derivatives are evaluated at 
$(x,0)$.  The Ricci curvature component $\Rt_{i\infty}$ of $\gt$ vanishes
to order $n/2-1$, and $D_i$ is given by:
$$
D_i = -2(\pa_\rho^{n/2-1}\Rt_{i\infty})|_{\rho =0,\, t=1}.   
$$
When $n=4$, this reduces to:
$$
D_i=4P^{jk}P_{ij,k}-3 P^{jk}P_{jk,i}+2P_i{}^jP^k{}_{k,j}. 
$$

If the metric $\gt$ in Theorem~\ref{main} is written in the form 
\eqref{formgt}, \eqref{normalvanishing} and $g_{ij}$ is also expanded as in  
\eqref{normalvanishing} (with expansion including a smooth term
$g_{ij}^{(0)}$), then  
$g^{(1)}_{ij}|_{\rho =0} =c_1 \cO_{ij}$ and 
$g_{i\infty}^{(1)}|_{\rho =0} = c_2 (A_{ij},{}^j -D_i)$ for 
nonzero constants $c_1$, $c_2$.    

Theorem~\ref{main} is proved by an inductive perturbation analysis of the
equation of vanishing Ricci curvature for metrics of the form given by 
Proposition~\ref{normalformprop}.  In the proof, one sees that the smooth
part $\gt^{(0)}$ of any inhomogeneous ambient metric is a smooth  
ambient metric for $(M,[g])$ in the sense of \S \ref{homo}. 

The ambiguity tensor has a well-defined transformation law under conformal
change.   
Let $\gt$ be the inhomogeneous ambient metric in normal form 
determined by initial metric $g$ and ambiguity tensor $A_{ij}$.   
Suppose we choose another metric $\wh{g}=e^{2\Up}g$ in the 
conformal class.  There is a uniquely determined smooth homogeneous
diffeomorphism $\Phi$ with $\Phi|_{\cG}=Id$ so that $\Phi^*\gt$ is in
normal 
form relative to $\wh{g}$.  Now $\Phi^*\gt$ determines an ambiguity  
tensor which we denote
by $\wh{A}_{ij}$, defined by the version of the equation \eqref{ambig} for 
the Taylor coefficient of the smooth part of $\Phi^*\gt$ in the
identification $\cGt = \R_+\times M\times \R$ induced by $\wh{g}$.    
We describe next the expression for $\wh{A}_{ij}$ in terms of 
$A_{ij}$ and $\Up$.  
This transformation law is best understood in terms of another 
trace-free symmetric 2-tensor on $M$ which is a modification of 
the ambiguity tensor.   

According to Proposition~\ref{normalformprop}, the smooth part 
$\gt_{IJ}^{(0)}$ of an inhomogeneous ambient metric $\gt_{IJ}$  
in normal form relative to $g$ takes the form 
$$
\gt_{IJ}^{(0)}=
\begin{pmatrix}
2\rho&0&t\\
0&t^2g_{ij}^{(0)}&0\\ 
t&0&0
\end{pmatrix},
$$
where $g_{ij}^{(0)}(x,\rho)$ is a smooth 1-parameter family of metrics on
$M$.  Since $\gt^{(0)}$ is itself a smooth ambient metric,
the derivatives   
$\pa_\rho^m g_{ij}^{(0)}|_{\rho =0}$ for $m< n/2$ and the trace 
$g^{kl}\pa_\rho^{n/2}g^{(0)}_{kl}|_{\rho =0}$ 
are the same  
natural tensors $\pa_\rho^m g_{ij}|_{\rho =0}$ and 
$g^{kl}\pa_\rho^{n/2}g_{kl}|_{\rho =0}$  
which occur in the expansion \eqref{expand}.  
Consider the value at    
$\rho =0$, $t=1$ of the component 
$\Rt^{(0)}_{\infty
  ij\infty,\,\underbrace{\scriptstyle{\infty\ldots\infty}}_{n/2-2}}$
of the iterated covariant derivative of the curvature tensor of  
$\gt_{IJ}^{(0)}$ .  This component depends on derivatives of orders $\leq
n/2$ of the components of $\gt_{IJ}^{(0)}$.  An inspection of the formula
for this covariant derivative component yields the following proposition. 

\begin{proposition}\label{curvderiv}
$$
\Rt^{(0)}_{\infty 
  ij\infty,\,\underbrace{\scriptstyle{\infty\ldots\infty}}_{n/2-2}}|_{\rho
  =0,\,t=1}  
=\tfrac12  (A_{ij} + K_{ij}),
$$
where $K_{ij}$ is a natural trace-free symmetric tensor which can be 
expressed algebraically in terms of
the tensors $\pa_\rho^m g_{ij}|_{\rho =0}$, $m<n/2$.  
\end{proposition}

The tensors $\pa_\rho^m g_{ij}|_{\rho =0}$, $m<n/2$, are determined by the 
inductive solution of the equation $\Ric(\gt)=0$ for a smooth ambient
metric in normal form, and 
$K_{ij}$ is expressed in terms of these via the formula for covariant
differentiation.  To the extent that $K_{ij}$, and therefore also  
its conformal transformation law, may be 
regarded as known, the conformal transformation law of $A_{ij}$ is
determined by that of $\Rt^{(0)}_{\infty 
  ij\infty,\,\underbrace{\scriptstyle{\infty\ldots\infty}}_{n/2-2}}|_{\rho
  =0,\,t=1}$.  Henceforth we shall write ${\bf A}_{ij}=A_{ij}+K_{ij}$ for  
this modified ambiguity tensor.  For $n=4$, one has 
$K_{ij}= -2\tf\left(P_i{}^kP_{jk}\right)$ so that 
\begin{equation}\label{Aform}
{\bf A}_{ij} = A_{ij} - 2\tf \left(P_i{}^kP_{jk}\right). 
\end{equation}

The conformal transformation law of ${\bf A}_{ij}$ can be expressed
succinctly.  Define the strength of lists of indices in $\R^{n+2}$ as
follows.  Set $\|0\| = 0$,  
$\|i\| = 1$ for $1\leq i \leq n$, and $\|\infty\|=2$.  For a list,  
write $\|I\ldots J\| = \|I\|+\cdots +\|J\|$.  An inductive analysis of the 
formula for covariant differentiation shows that if $r\geq 0$ and 
$\|ABCDF_1\ldots F_r\|\leq n+1$, then the curvature component 
$\Rt^{(0)}_{ABCD,F_1\cdots F_r}|_{\rho =0,\,t=1}$  
defines a natural tensor on $M$  
as the indices between $1$ and $n$ vary and those which are $0$ or $\infty$
remain fixed.  Set
$$
p^A{}_I = 
\begin{pmatrix}
1&\Up_i&-\frac12 \Up_k\Up^k\\
0&\delta^a{}_i&-\Up^a\\
0&0&1
\end{pmatrix}.
$$

\begin{proposition}\label{Atrans}
Set $r=n/2-2$.  
Under the conformal change $\wh{g}=e^{2\Up}g$, the modified ambiguity
tensor transforms by:
$$
e^{(n-2)\Up}\wh{\bf A}_{ij} = {\bf A}_{ij} + 
2 \sum{\vphantom{\sum}}'\Rt^{(0)}_{ABCD,F_1\cdots F_r}|_{\rho =0,\,  
  t=1}  
p^A{}_\infty p^B{}_ip^C{}_jp^D{}_\infty p^{F_1}{}_\infty\cdots
p^{F_r}{}_\infty,
$$
where $\sum'$ indicates the sum over all indices except 
$ABCDF_1\cdots F_r =\infty ij \infty \infty\cdots \infty$.  
\end{proposition}

The upper-triangular form of $p^A{}_I$ implies that 
any component $\Rt^{(0)}_{ABCD,F_1\cdots F_r}|_{\rho =0,\, t=1}$ 
which occurs with nonzero coefficient in $\sum'$ necessarily satisfies  
$\|ABCDF_1\cdots F_r\|$ $\leq n+1$.  So all these components are natural 
tensors, which can be calculated algorithmically using the expansion 
of $g_{ij}(x,\rho)$ through order $< n/2$.  Thus Proposition~\ref{Atrans}
expresses the conformal 
transformation law of ${\bf A}_{ij}$ in terms of ``known'' natural tensors
and $\Up$ and its first derivatives.  For example, for $n=4$ this becomes 
$$
e^{2\Up}\wh{\bf A}_{ij} ={\bf A}_{ij} - 2\Up^l(C_{ijl}+C_{jil})
+2\Up^k\Up^l W_{kijl}.  
$$

Fix an even integer $n\geq 4$.  In dimension $d>n$, there is a trace-free
symmetric natural tensor satisfying the  
transformation law of Proposition~\ref{Atrans}; namely 
$2\Rt^{(0)}_{\infty 
  ij\infty,\,\underbrace{\scriptstyle{\infty\ldots\infty}}_{n/2-2}}|_{\rho
  =0,\,t=1}$.  For instance,  
$$
\Rt^{(0)}_{\infty ij\infty}|_{\rho =0,\,t=1}=-\frac{1}{d-4}B_{ij}
$$   
if $d\neq 4$.  Considered formally in the dimension $d$, the component 
$\Rt^{(0)}_{\infty 
  ij\infty,\,\underbrace{\scriptstyle{\infty\ldots\infty}}_{n/2-2}}|_{\rho
  =0,\,t=1}$
has a  
simple pole at $d=n$ whose residue is a multiple of $\cO_{ij}$. 
The formal continuation to $d=n$ of the transformation law for this
component 
is the statement of conformal invariance of the obstruction tensor.  The   
modified ambiguity tensor ${\bf A}_{ij}$ provides a substitute for the
natural tensor 
$\Rt^{(0)}_{\infty 
  ij\infty,\,\underbrace{\scriptstyle{\infty\ldots\infty}}_{n/2-2}}|_{\rho
  =0,\,t=1}$, 
which does not occur in dimension $n$ because of 
the conformal invariance of the obstruction tensor.  

The transformation law of Proposition~\ref{Atrans} can be reinterpreted in
terms of tractors.  Recall (see, for example, \cite{BEGo}) that a cotractor
field of weight $w$ on a conformal manifold $(M,[g])$ can be expressed upon
choosing a 
representative metric $g$ as $v_I = (v_0,v_i,v_\infty)$, where $v_0$,
$v_\infty$ are functions on $M$ and $v_i$ is a 1-form, such that under the 
conformal change $\wh{g}=e^{2\Up}g$, $v_I$ transforms by
\begin{equation}\label{tractor}
\wh{v}_I = e^{(w+1-2\delta^\infty_I)\Up}v_Ap^A{}_I. 
\end{equation}
For $r<n/2-2$, the components of the tensor
$\tilde{\nabla}^r \Rt^{(0)}|_{\rho =0,\, t=1}$ define natural tensors on
$M$ and satisfy the correct transformation laws under conformal change 
to  
define a (natural) cotractor field of rank $4+r$ and weight $-2-r$.  
(A detailed discussion of the relation between ambient    
curvature and tractors can be found in \cite{CG}.) 
For $r=n/2-2$, the same is true of all components except for   
$\Rt^{(0)}_{\infty 
  ij\infty,\,\underbrace{\scriptstyle{\infty\ldots\infty}}_{n/2-2}}|_{\rho
  =0,\,t=1}$ (and components obtained from this one via the symmetries of
the curvature tensor).  The transformation law of  
Proposition~\ref{Atrans} is precisely that required by \eqref{tractor} so
that the prescription  
$\Rt^{(0)}_{\infty 
  ij\infty,\,\underbrace{\scriptstyle{\infty\ldots\infty}}_{n/2-2}}|_{\rho
  =0,\,t=1}= \tfrac12 {\bf A}_{ij}$ completes 
$\tilde{\nabla}^{n/2-2} \Rt^{(0)}|_{\rho =0,\, t=1}$ to a cotractor field
of rank $n/2+2$ and weight $-n/2$.  
Thus, choosing an ambiguity tensor $A_{ij}$ is entirely equivalent to   
completeing $\tilde{\nabla}^{n/2-2} \Rt^{(0)}|_{\rho =0,\,t=1}$ to a
cotractor field.    

When $n$ is even, one can construct homogeneous ambient metrics with
asymptotic expansions involving $\log{|r|}$ (\cite{K}).  Such nonsmooth  
homogeneous ambient metrics all have asymptotic expansions of the form 
$\sum_{N\geq 0} \bar{g}^{(N)} (r^{n/2}\log |r|)^N $, where  
the 
$\bar{g}^N$ are smooth and homogeneous of degree 2 ([FG2]).  This
suggests that the inhomogeneous ambient metrics considered here might
actually have expansions of the form 
\begin{equation}\label{modifiedgt}
\sum_{N\geq 0} \bar{g}^{(N)} (r^{n/2}\log |r_\#|)^N, 
\end{equation}
again with smoooth homogeneous coefficients $\bar{g}^{(N)}$.  
But this is not the case: for inhomogeneous ambient metrics, the
coefficient 
$\gt^{(2)}$ does not in general vanish at $\rho =0$.  In fact, we have
\begin{proposition}\label{log2}
Let $\gt_{IJ}$ be an inhomogeneous ambient metric in normal form.  When
$\gt_{IJ}$ is written in the form \eqref{formgt} and its components are
expanded as in \eqref{normalvanishing} (with a corresponding expansion for
$g_{ij}$ including a smooth term
$g_{ij}^{(0)}$), then $g_{ij}^{(2)}|_{\rho =0}=0$,  
$g_{i\infty}^{(2)}|_{\rho =0}=0$, and 
$g_{\infty\infty}^{(2)}|_{\rho =0} = c\,\cO_{ij}\cO^{ij}$, for a 
nonzero constant $c$.  Thus if $\gt$ has the form  
\eqref{modifiedgt}, then $\cO^{ij}\cO_{ij}=0$.  
Conversely, if $\cO^{ij}\cO_{ij}=0$, then $\gt$ has the form
\eqref{modifiedgt} for any choice of ambiguity.
\end{proposition}

In Proposition~\ref{log2}, one would expect that there is a condition
analogous to  
$\cO^{ij}\cO_{ij}=0$ arising from the coefficient of each of $(\log |2\rho 
t^2|)^N$  for $N\geq 2$.  The fact that all of these are satisfied once the
condition holds for $N=2$ is surprising.    

Proposition~\ref{log2} raises the question of whether there are interesting
classes of conformal manifolds for which $\cO_{ij}$ is nonvanishing but 
$\cO^{ij}\cO_{ij}=0$.  Of course, this cannot happen in definite
signature.  But the Fefferman conformal structure of any nondegenerate 
integrable CR manifold satisfies $\cO^{ij}\cO_{ij}=0$, and typically they
satisfy also that $\cO_{ij}$ is nonvanishing.  These observations follow
from an identification of the obstruction tensor of a Fefferman metric.   

The analogue of the obstruction tensor in CR geometry is the scalar CR
invariant, denoted here by $L$, which obstructs the existence of 
smooth solutions of Fefferman's complex 
Monge-Amp\`ere equation.   It is defined as follows.  Let 
$\cM\subset \C^n$ be a hypersurface with smooth defining function $u$ whose
Levi form $-u_{i\bar{j}}|_{T^{1,0}\cM}$ has signature $(p,q)$,  $p+q=n-1$.    
Fefferman showed in \cite{F1} that there is such a $u$  
uniquely determined mod $O(u^{n+2})$ such that $J(u)=1+O(u^{n+1})$, where 
$$
J(u)= (-1)^{p+1}\det
\begin{pmatrix}
u&u_{\bar{j}}\\
u_i&u_{i\bar{j}}
\end{pmatrix}_{1\leq i,j \leq n}.
$$
The invariant $L$ is defined to be a constant multiple of
$(J(u)-1)/u^{n+1}|_\cM$, and is 
independent of the choice of $u$ satisfying $J(u)=1+O(u^{n+1})$.    
The fundamental properties of $L$ are derived in   
\cite{L}, \cite{Gr1}, \cite{Gr2}. 

\begin{proposition}\label{cr}
Let $\theta$ be a pseudohermitian form for an integrable nondegenerate CR
manifold $\cM$ and let 
$g$ be the associated representative of the Fefferman conformal structure
on the circle bundle $\cC$.  Then the obstruction tensor
$\cO$ of 
$g$ is a nonzero constant multiple of the pullback to $\cC$ of $L\theta^2$.    
\end{proposition}

\begin{proof}
It suffices to assume that $\cM\subset \C^n$ is  
embedded.  The circle bundle $\cC = S^1\times \cM$ has dimension $2n$.  Let
$u$ be a smooth defining function for $\cM$ satisfying 
$J(u)=1+O(u^{n+1})$ as above.  According to Fefferman's original
definition, the representative $g$ associated to 
$\theta = \frac{i}{2}(\pa u-\bar{\pa} u )$ is the pullback to  
$$
S^1\times
\cM = \{(z^0,z): |z^0|=1,\,z\in \cM\}\subset \C^*\times \C^n
$$ 
of the  K\"ahler metric $\gt$ on a neighborhood of $\C^*\times \cM$ in
$\C^*\times \C^n$  given by  
\begin{equation}\label{crambient}
\gt=\pa^2_{I{\overline J}}(-|z^0|^2u)dz^Id{\overline z}^J.
\end{equation}
The Ricci curvature of $\gt$ is 
$$
\pa^2_{I{\overline J}}(\log J(u))dz^Id{\overline z}^J 
= c L u^{n-1}\pa u {\overline \pa}u+O(u^n).
$$
Since the pullback to $\cG=\C^*\times \cM$ of $\pa u {\overline \pa}u$ is a
multiple 
of $\theta^2$, which is trace-free with respect to $g$, it is evident that 
this Ricci curvature is $O^+_{IJ} (\rho^{n-1})$.   
Now $\gt$ is clearly homogeneous, and it satisfies the initial condition 
\eqref{initial} by the definition of the Fefferman metric $g$.  Therefore 
$\gt$ is a smooth   
ambient metric for the Fefferman conformal structure in the sense of our
earlier definition.  Now $T=2\Re z^0\pa_{z^0}$, and an easy calculation
gives  
$\gt(T,T) = -|z^0|^2u$.  The definition \eqref{obsdef} of the obstruction 
tensor thus reduces to a constant multiple of $L\theta^2$.          
\end{proof}

Since $\theta$ is null with respect to $g$, 
Proposition~\ref{cr}  implies 
that the obstruction tensor of a Fefferman metric satisfies
$\cO_i{}^k\cO_{kj}=0$.  In particular $\cO^{ij}\cO_{ij}=0$, so Fefferman
metrics satisfy the condition in Proposition~\ref{log2}.  

Upon applying Theorem~\ref{main} to the Fefferman conformal structure of a
CR manifold, one obtains a
family of inhomogeneous ambient metrics with ambiguity an 
arbitrary trace-free symmetric tensor on the circle bundle.  In \cite{H}, a
family of ``inhomogeneous CR ambient metrics'' associated to a
nondegenerate 
hypersurface in $\C^n$ was constructed with ambiguity a scalar function on
the hypersurface.  The relation between these constructions is not clear.  
The metrics constructed in \cite{H} are not straight in general.  It seems  
reasonable to guess that there are inhomogeneous, nonsmooth
diffeomorphisms with expansions involving $\log |r_\#|$ relating the
constructions (for an appropriately restricted ambiguity in the conformal
construction).  We intend to investigate this issue in the future.

%

A homogeneous ambient metric for a conformal manifold $(M,[g])$ gives rise
to an asymptotically hyperbolic ``Poincar\'e'' metric by a procedure
described in \cite{FG1}.  Namely, if $\gt$ is a homogeneous ambient metric,
then the pullback $g_+$ of $\gt$ to the hypersurface $\{\gt(T,T)=-1\}$ is a  
metric of signature of $(p+1,q)$ which in suitable coordinates is
asymptotically hyperbolic with conformal infinity $(M,[g])$.  The
homogeneity of $\gt$ implies that the condition
$\Ric(\gt)=0$ is equivalent to $\Ric(g_+)=-ng_+$.  The same construction
can also be carried out for an inhomogeneous ambient metric. 
Proposition~\ref{Tsmooth} shows that $\gt(T,T)$ is still a smooth
homogeneous defining function.  When we take $r_\#=\gt(T,T)$,
all log terms in the asymptotic expansion \eqref{metricexpansion} 
vanish on the hypersurface $\gt(T,T)=-1$, so the Poincar\'e metric $g_+$
which 
one obtains agrees with that defined by the smooth part $\gt^{(0)}$.
In particular, $g_+$ has a smooth conformal compactification with no log
terms.  The $g_+$ which arise as $\gt$ varies over the family of
inhomogeneous ambient metrics 
associated to $(M,[g])$ form a family of smoothly compactifiable
Poincar\'e metrics invariantly associated to $(M,[g])$ up to choice of
ambiguity and up to a smooth  
diffeomorphism restricting to the identity on $M$.  However, 
$\gt^{(0)}$ is not usually Ricci flat, so $g_+$ will not 
in general be Einstein.  We do not know a direct characterization of either 
the smooth parts $\gt^{(0)}$ or the resulting family of Poincar\'e metrics
$g_+$, other than to  
say that they arise from inhomogeneous ambient metrics by taking smooth 
parts.

\section{Scalar invariants}\label{invariants}

An original motivation for the ambient metric construction in \cite{FG1}
was to construct scalar conformal invariants.  In even dimensions, the
construction of invariants in \cite{FG1} terminates at finite order owing
to the failure 
of the existence of infinite order smooth ambient metrics.  Inhomogeneous
ambient metrics can be used to extend the construction of invariants to  
all orders.  The new invariants generally depend both on the 
initial metric and the ambiguity.

By a scalar invariant $I(g)$ of metrics we mean a polynomial in 
the variables $(\partial^{\alpha}g_{ij})_{|\alpha| \geq 0}$ 
and $|\det{g_{ij}}|^{-1/2}$, which is coordinate-free in the sense that its  
value is independent of orientation-preserving changes of the coordinates
used to express and differentiate $g$.  Such a scalar invariant 
is said to be {\it even} if it is also unchanged under
orientation-reversing changes of coordinates, and {\it odd} if it changes
sign under orientation-reversing coordinate changes.  
It is said to be conformally invariant 
of weight $w$ if $I(\Omega^2 g) = \Omega^{w}I(g)$ for smooth positive
functions $\Omega$.  By a scalar conformal invariant we will mean a
scalar invariant of metrics which is 
conformally invariant of weight $w$ for some $w$.  

First recall the characterization of scalar conformal invariants in odd
dimensions.  Denote by $\nt^r\Rt$ the $r$-th iterated covariant derivative
of the curvature tensor of a smooth ambient metric, by $\mut\in
\bigwedge^{n+2}T^*\cGt$ the volume 
form of $\gt$, and set $\mut_0 = T\into \mut$.  Consider complete
contractions of these tensors:
\begin{gather}\label{ambcontr}
\begin{gathered}
\contr(\nt^{r_1}\Rt\otimes \cdots \otimes \nt^{r_L}\Rt)\\
\contr(\mut\otimes \nt^{r_1}\Rt\otimes \cdots \otimes \nt^{r_L}\Rt)\\ 
\contr(\mut_0\otimes \nt^{r_1}\Rt\otimes \cdots \otimes \nt^{r_L}\Rt), 
\end{gathered}
\end{gather}
where the contractions are taken with respect to $\gt$.  
Each such complete contraction defines a homogeneous function on
$\cGt$ whose restriction to $\cG$  is independent of the diffeomorphism
indeterminacy in the smooth ambient metric.    
If $g$ is a metric in the conformal class, evaluating this homogeneous
function at the image of $g$ viewed as a section of $\cG$
defines a function on $M$ which depends polynomially on the derivatives of
$g$. Every such function is a scalar conformal invariant.   
Contractions of the first type in \eqref{ambcontr} define even invariants 
and contractions of 
the second and third types define odd invariants.  
A linear combination of such complete contractions, all having the 
same weight, is called a Weyl invariant.
The characterization of conformal invariants for $n$ odd states that  
every scalar conformal invariant is a Weyl invariant.
This follows from the invariant theory of \cite{BEGr} together with  
the jet isomorphism theorem of \cite{FG2}, as outlined in \cite{BEGr}. 
Full details will appear in \cite{FG2}. 

When $n$ is even, we carry out the same construction, but now
replacing the $\nt^r \Rt$ in the complete contractions \eqref{ambcontr} by
the corresponding covariant derivatives of 
curvature of the smooth part of an inhomogeneous ambient metric.  
Of course, the resulting functions on $M$ generally depend on the choice of
ambiguity  
tensor as well as the metric $g$; they can be regarded as conformal
invariants of the pair $(g,A)$.  We still refer to linear combinations of
such complete contractions as Weyl invariants.  Some such Weyl invariants 
may actually be independent of the ambiguity, in which case they define
scalar conformal invariants.  We call such Weyl invariants
ambiguity-independent Weyl invariants.  We do not have a systematic general
way of constructing ambiguity-independent Weyl invariants, but we do have
examples and conditions for Weyl invariants to be ambiguity-independent in
interesting cases.  Also, we can extend the jet isomorphism theorem of
\cite{FG2} and the invariant theory of \cite{BEGr} to prove 
that all scalar conformal invariants are of this form if $n\equiv 2 \mod
4$:   

\begin{theorem}
If $n\equiv 2 \mod 4$, then every scalar conformal invariant is an
ambigu-
\newline
ity-independent Weyl invariant.  
\end{theorem}

If $n$ is a multiple of 4, there are exceptional odd invariants which do 
not arise as ambiguity-independent Weyl invariants.  The existence of 
exceptional invariants for a linearized problem was observed in    
\cite{BEGr} and in this context they were studied systematically in  
\cite{BG}.  We follow the approach of \cite{BG} to complete the description
of scalar conformal invariants when $n\equiv 0 \mod 4$ as follows.  

Choose a representative metric $g$ and corresponding identification 
$\cGt = \R_+\times M\times \R$.  Define an $n$-form $\eta$ on 
$\cGt$ by $\eta = t^{-2}\pa_\rho \into \mu_0$.    
Consider a partial contraction
$$
L=\text{pcontr}(
\underset{n/2}{\eta\otimes\underbrace{\Rt\otimes \cdots \otimes \Rt}}), 
$$
where all the indices of $\eta$ are contracted.
While $\eta$ depends on the choice of $g$, it can be shown that $L|_\cG$  
is independent of this choice and also of the choice of 
ambiguity tensor (for $n>4$ this is obvious since 
$\Rt|_\cG$ is already independent of the ambiguity).    
A  scalar invariant can be obtained from $L$ by applying an ambient 
realization of the 
tractor $D$ operator (see again \cite{CG}).  This is a differential  
operator $D$ defined initially as a map  
$D:\Gamma(\otimes^pT^*\cGt)(w)\rightarrow \Gamma(\otimes^{p+1}T^*\cGt)(w)$,
where $\Gamma(\otimes^pT^*\cGt)(w)$ denotes the space of covariant
$p$-tensor fields homogeneous of degree $w$ in the sense of the
previous sections ($\delta^*_s f = s^wf$), by:
$$
Df=\nt f-\frac{1}{2(n+2(w-p-1))}dr_\#\otimes\widetilde{\Delta} f.
$$
Here $\nt$ and $\widetilde{\Delta}=\nt^I\nt_I$ are defined with respect to 
$\gt^{(0)}$ and it is assumed that $w\neq p+1-n/2$.  One checks that $D$
acts tangentially to $\cG$ in the sense 
that $Df$ vanishes on $\cG$ if $f$ does, so $D$ induces a map
$$
D:\Gamma(\otimes^pT^*\cGt|_{\cG})(w)\rightarrow
\Gamma(\otimes^{p+1}T^*\cGt|_{\cG})(w). 
$$  
The homogeneities are such that the expression
$$
D^I D^J\cdots D^KL_{(IJ\cdots K)} 
$$
is defined and the above discussion implies that it restricts to 
$\cG$ to give a homogeneous function independent of the choice of ambiguity
tensor.  (Here $(IJ\cdots K)$ denotes symmetrization over the enclosed
indices.)  As in the case of Weyl invariants, when
evaluated at the image of $g$ as a section of $\cG$, one obtains 
a scalar conformal invariant.  
Such conformal invariants are odd and are called  {\it basic exceptional 
invariants}.  It is easy to see that there are only finitely many basic 
exceptional invariants in 
any dimension $n\equiv 0 \mod 4$ (the construction works just as well when
$n\equiv 2\mod 4$, but all basic exceptional invariants vanish in this
case).  For $n=4$ there are only two nonzero basic exceptional invariants:
$$
\eta^{IJKL}\Rt_{IJ}{}^{AB}\Rt_{KLAB}
$$ 
of weight $-4$, and 
$$
D^AD^B(\eta^{IJKL}\Rt_{IJAC}\Rt_{KLB}{}^C)
$$ 
of weight $-6$.    
The first is a multiple of $|W^+|^2-|W^-|^2$, where $W^\pm$ denote the 
$\pm$ self-dual parts of the Weyl tensor, and the second is due to Bailey,
Eastwood, and Gover (see p. 1207 of \cite{BEGo}).  

By extending the jet isomorphism theorem and invariant theory of
\cite{FG2}, \cite{BEGr} and \cite{BG}, we can prove:    
\begin{theorem} 
If $n\equiv 0 \mod 4$,  then every even scalar conformal invariant is an
ambiguity-independent Weyl invariant, and every odd scalar conformal
invariant is a linear combination of an ambiguity-independent Weyl
invariant and basic exceptional invariants.    
\end{theorem}

The  above constructions use only the smooth part $\gt^{(0)}$ of an
inhomogeneous ambient metric.  It is also possible to construct 
invariants using the tensors $r^{(n/2-1)N+1}\gt^{(N)}$ for $N\geq 1$ in
\eqref{metricexpansion} (with $r_\#=\gt(T,T)$).  However, the above
theorems show that this is not necessary:  the invariants already can be
constructed using just the smooth part.  

We conclude with some examples.  The smooth part $\gt^{(0)}$ of an  
inhomogeneous ambient metric is not Ricci flat; the leading part of its 
Ricci tensor can be identified with the obstruction tensor.  So 
Weyl invariants involving the Ricci tensor of $\gt^{(0)}$ give rise to 
conformal invariants involving the obstruction tensor.  
Define an ambient version of the obstruction tensor by   
$$
\cOt_{IJ} = \tfrac{1}{n-2}\widetilde{\Delta}^{n/2-1}\Rt_{IJ},
$$
where $\tilde{\Delta}$ and $\Rt_{IJ}$ are with respect to $\gt^{(0)}$.   
One checks that $T\into \cOt = 0$ on $\cG$ and that $\cOt|_{T\cG}$ is the
tensor on $\cG$ homogeneous of degree $2-n$ defined by $\cO$.  Otherwise
put, in the identification $\cGt = \R_+\times M\times \R$ induced by a
representative $g$, one has $\cOt_{0I}=0$ and $\cOt_{ij}=t^{2-n}\cO_{ij}$ 
at $\rho =0$.  {From} this it is easy to see that the contraction  
$\cOt_{IJ}\cOt^{IJ}$ determines the conformal invariant $\cO_{ij}\cO^{ij}$, 
and similarly that $\Rt^{IJKA}\Rt_{IJK}{}^B\cOt_{AB}$ determines 
the conformal invariant $W^{ijka}W_{ijk}{}^b\cO_{ab}$.  A more interesting
example is to consider $\Rt^{IJKL}\cOt_{IK,JL}$.  One finds that for even
$n\geq 6$, the restriction of this quantity to $\cG$ is
ambiguity-independent and determines the conformal invariant   
$$
W^{ijkl}\cO_{ik,jl} -(n-1)W^{ijkl}P_{ik}\cO_{jl} +2n C^{jkl}\cO_{jk,l} 
+\tfrac{n(n-1)}{n-4} B^{jk}\cO_{jk}.
$$
For $n=4$, the Weyl invariant determined by $\Rt^{IJKL}\cOt_{IK,JL}$
depends on the ambiguity tensor.  It is given by the same formula,
except that in the last term $\frac{1}{n-4}B^{jk}$ is replaced by 
$-\frac12 {\bf A}^{jk}$ with ${\bf A}$ given by \eqref{Aform}.

\end{document}